\documentclass[12pt, reqno]{amsart}
\usepackage{inputenc}
\usepackage{amsmath, amsthm, amscd, amsfonts, amssymb, graphicx, color, hyperref}
\newtheorem{theorem}{Theorem}[section]

\newtheorem{proposition}[theorem]{Proposition}
\newtheorem{corollary}[theorem]{Corollary}
\theoremstyle{definition}
\newtheorem{definition}[theorem]{Definition}
\newtheorem{example}[theorem]{Example}

\theoremstyle{remark}
\newtheorem{remark}[theorem]{Remark}
\numberwithin{equation}{section}
\begin{document}
\setcounter{page}{1}

\title[ Left and right generalized Drazin invertible operators ]{Left and right generalized Drazin invertible operators and Local spectral theory }

\author[M. Benharrat, K. Miloud Hocine and B. Messirdi]{ Mohammed Benharrat$^1$  Kouider Miloud Hocine$^2$,    and Bekkai Messirdi$^3$
}

\address{$^{1}$ Laboratoire de Math\'ematiques Fondamentales et Appliqu\'ees d'Oran (LMFAO),  D\'epartement de Math\'ematiques et Informatique, 
Ecole National Polythecnique d'Oran (Ex. ENSET d'Oran);  B.P. 1523  Oran-El M'Naouar, 31000 Oran, Alg\'erie.}

\email{\textcolor[rgb]{0.00,0.00,0.84}{ mohammed.benharrat@enp-oran.dz}}

\address{$^{2}$,$^3$ Laboratoire de Math\'ematiques Fondamentales et Appliqu\'ees d'Oran (LMFAO),  D\'epartement de Math\'ematiques , Universit\'e d'Oran 1, 31000 Oran, Alg\'erie.}
\email{\textcolor[rgb]{0.00,0.00,0.84}{miloud.houcine@edu.univ-oran1.dz\\
messirdi.bekkai@univ-oran.dz}}




\begin{abstract}
In this paper, we give some characterizations of the left and right generalized Drazin invertible bounded operators in Banach spaces by means of the single-valued extension property (SVEP). In particular, we show that a bounded operator is left (resp. right) generalized Drazin invertible if and only if admits a generalized Kato decomposition and has the SVEP at 0 (resp. it admits a generalized Kato decomposition  and its adjoint has the SVEP at 0. In addition,  we prove that both of  the left and the  right generalized Drazin    operators are  invariant under  additive commuting  finite rank  perturbations. Furthermore, we investigate the transmission of some local spectral properties from a bounded linear operator, as the SVEP, Dunford property $(C)$, and property $(\beta)$, to its generalized Drazin inverse.

\end{abstract} 
\maketitle
\section{Introduction}
The generalized Drazin inverse for operators  arises naturally in the context of isolated
spectral points and  becomes  a theoretical and practical tool in algebra and analysis (Markov chains, singular differential and difference equations, iterative methods...). 
The Drazin inverse was originally defined in 1958  for semigroups (\cite{Drazin}).  When $\mathcal{L}(X)$ is  the Banach algebra of all bounded linear operators acting on an infinite-dimensional  complex Banach space $X$, then $S\in \mathcal{L}(X)$ is the Drazin inverse of $T\in \mathcal{L}(X)$ if 
\begin{equation}\label{defDrazin}
ST=TS \quad STS=S \hbox{ and } TST=T+U \text{ where $U$ is a nilpotent operator}.
\end{equation}

The concept of Drazin invertible  operators has been generalized by Koliha (\cite{Koliha0}) by replacing  the nilpotent operator $U$ in \eqref{defDrazin} by a quasinilpotent operator. In this case, $S$ is called a generalized Drazin inverse of $T$. Note that this extension was anticipated by Harte in \cite{Harte}. Recently, in \cite{miloud}, the authors introduced the left and the right generalized Drazin invertible operators. These two classes of operators are a continuation and refinement of the research treatment of the Drazin inverse  in  Banach space operators theory. It proved that an operator $T\in \mathcal{L}(X)$ is   left (resp.  right) generalized Drazin inverse if and only if $T=T_1 \oplus T_2$ where $T_1$ is bounded below (resp. surjective) and $T_2$ is quasinilpotent operator. Furthermore, these operators are characterized via the isolated points of the approximate point spectrum (resp. surjective spectrum)   \cite[ Theorem 3.8; Theorem 3.10]{miloud}.

The main objective of this paper is to continue studying these operators via the local spectral theory. In Section \ref{sec:Prel}, we give some preliminary results which our investigation will be need. In Section \ref{sec:3}, we present many new and interesting characterizations of  the left  (resp. the right) generalized Drazin invertible  operators  in terms of the generalized Kato decomposition and  the single-valued extension property. We also show that an operator admits a generalized Kato decomposition and has the SVEP at 0 is precisely left generalized Drazin invertible and conversely. Similarly, an operator $T$ is right generalized Drazin invertible if and only if $T$ admits a generalized Kato decomposition and its adjoint $T^*$ has the SVEP at 0. In particular, we prove that the left generalized Drazin spectrum and right generalized Drazin spectrum of a bounded operator are invariant under   commuting  finite rank  perturbations. In section \ref{sec:4},  we study the relationships between the local spectral properties of an operator  and
the local spectral properties of its generalized Drazin inverse, if this exists. In particular, a reciprocal relationship analogous to spectrum of invertible operator and its inverse, is established between the nonzero points  of the local spectrum of a generalized Drazin invertible operator having SVEP and the nonzero points of the local spectrum of its generalized  Drazin inverse. We also show that many local
spectral properties, as SVEP, Dunford property $(C)$, property $(\beta)$, property $(Q)$ and  decomposability, are transferred from a generalized Drazin invertible operator  to its generalized Drazin inverse. This section extends the results of \cite{AienaT} from the case of Drazin invertible operators to case of the generalized Drazin invertible operators. Finally, by a counterexample we show that these local
spectral properties are  not transferred in the case of the left  (resp. the right) generalized Drazin invertible  operators.

\section{Preliminaries} \label{sec:Prel}
 Let  $\mathcal{L}(X)$ be  the Banach algebra of all bounded linear operators acting on an infinite-dimensional  complex Banach space $X$. 
 For  $T\in \mathcal{L}(X)$ write $ N(T)$, $ R(T)$,  $\sigma(T)$  and $ \rho(T)$ respectively, the null space, the range, the spectrum and the resolvent set of $T$. The nullity and the deficiency of $T$ are defined respectively by $\alpha(T) = \text{dim} N(T)$ and $\beta(T) = \text{dim} X/R(T).$ Here $I$ denote the identity operator in $X$. By iso$\sigma(T)$ and acc$\sigma (T)$ we define the set of all isolated and accumulation spectral points of $T$.

If $M$ is a subspace of $X$ then  $T_{M}$ denote the restriction of $T$ in $M$. Assume that $M$ and $N$ are two closed subspaces of $X$ such that $X=M\oplus N$ (that is $H=M+N$ and $M\cap N={0}$). We say that $T$ is completely reduced by the pair $(M,N)$, denoted as $(M, N) \in Red(T)$, if $T(M)\subset M$, $T(N)\subset N$ and $T=T_{M}\oplus T_{N}$. In such case we have   $N(T)=N(T_{M})\oplus N(T_{N})$, $R(T)=R(T_{M})\oplus R(T_{N})$,  and $T^{n}=T_{M}^{n}\oplus T_{N}^{n}$ for all $n\in \mathbb{N}$. An operator is said to be bounded below if it is injective with closed range.
 
Recall that (see, e.g. \cite{Kaashoek}) the ascent $a(T)$ of an operator $T\in\mathcal L(X)$ is defined as the smallest nonnegative integer $p$ such that $N(T^p)=N(T^{p+1})$. If no such an integer exists, we set $a(T)=\infty$. Analogously, the smallest nonnegative integer $q$ such that $R(T^q)=R(T^{q+1})$ is called the descent of $T$ and denoted by $d(T)$. We set $d(T)=\infty$ if for each $q$, $R(T^{q+1})$ is a proper subspace of $R(T^q)$. It is well known that if the ascent and the descent of an operator are finite, then they are equal.  
 
 Associated with an operator $T\in \mathcal{L}(X)$  there are  two (not necessarily closed) linear subspaces of $X$ invariant by $T$, played an important role in the development of  the generalized  Drazin inverse for $T$, the quasinilpotent part $H_{0}(T )$ of $T$:

$$H_{0}(T ) =\left\{ x \in x: \lim_{n\rightarrow \infty} \left\|T^{n}x\right\|^{\frac{1}{n}}=0 \right\},  $$
and the analytical core $K(T )$ of $T$:
$$ K(T ) =\{ x \in X: \text{ there exist a sequence } (x_n) \text{ in } X \text{ and a constant } \delta >0 $$ 
$$ \text{ such that } T x_{1} = x, T x_{n+1} = x_{n} \text{ and } \|x_n\| \leq \delta^{n} \|x\| \text{ for all } n \in \mathbb{N} \}.$$
It is well-known that  necessary and sufficient condition for $T\in \mathcal{L}(X)$ to be generalized Drazin invertible is that 0$\notin$acc$\sigma (T)$. Equivalently,
$K(T )$ and $H_{0}(T) $ are both closed, $X = H_{0}(T )\oplus K(T)$, the restriction of $T$ to  $H_{0}(T )$ is a quasinilpotent operator, and the restriction of $T$ to $K(T)$ is invertible,

Recently, by the use  of this two subspaces, in \cite{miloud}, the authors defined and studied a new class of operators called left and right generalized Drazin invertible operators as a generalization of left and right Drazin invertible operators.
\begin{definition}\label{def 1.1}  An operator $T\in\mathcal L(X)$ is said to be  right generalized Drazin invertible if $K(T)$ is closed and complemented with a subspace $N$ in $X$ such that $T(N)\subset N \subseteq H_{0}(T)$.
\end{definition}
\begin{definition}\label{def 1.2} An operator $T\in\mathcal L(X)$ is said to be  left generalized Drazin invertible if $H_{0}(T)$ is closed and complemented with a subspace $M$ in $X$ such that $T(M)\subset M$ and $T(M)$ is closed.
\end{definition}
 Examples of left generalized Drazin  invertible operators are the operators of the following classes:
\begin{itemize}
	\item Left  invertible operators.
	\item  Left Drazin invertible operators,
	\[LD(X)=\{T\in\mathcal L(X): a(T) \text{ is finite and }R(T^{a(T)+1})\text{ is closed}\}.\]
	In this case, $H_0(T)=N(T^p)$ with $a(T)=p<\infty$.
	\item Drazin invertible operators.
	\item  A   bounded paranormal operator $T$ on Hilbet space $X$ such that  $R(T)+H_{0}(T)$ and $H_{0}(T)$ are closed (see \cite[Proposition 3.15]{miloud}).
\end{itemize}
Fxamples of right generalized Drazin  invertible operators are the operators of the following classes:
\begin{itemize}
	\item  Right invertible  operators.
	\item Right Drazin invertible operators,
	 \[RD(X)=\{T\in\mathcal L(X): d(T) \text{ is finite and }R(T^{d(T)})\text{ is closed}\}.\]
	 In this case  $K(T)=R(T^q)$, with $d(T)=q<\infty$.
	\item Drazin invertible operators, 
	$$LD(X)\cap RD(X).$$
\end{itemize}
According to the Definitions \ref{def 1.1} and \ref{def 1.2}, we also have 

 Invertible operator  $\Longrightarrow$  Generalized Drazin invertible operator $\Longrightarrow$ Right (resp. Left) generelazed Drazin invertible operator.
 
 In the sequel the terms left (resp. right) generalized Drazin  invertible operator is used for the nontrivial case of the bounded below (resp. surjective) operators.

 The left Drazin spectrum, the right Drazin spectrum, the  Drazin spectrum, the  generalized Drazin spectrum, the left generalized Drazin spectrum and  the right generalized Drazin spectrum  of $T$ are, respectively, defined by\\

$\sigma_{lD}(T):=\{\lambda\in\mathbb{C}: \lambda I-T \notin LD(X)\},$\\

$\sigma_{rD}(T):=\{\lambda\in\mathbb{C}: \lambda I-T \notin RD(X)\},$\\

$\sigma_D(T)=\{\lambda\in\mathbb{C}:  \lambda I-T  \notin LD(X)\cap RD(X) \},$\\

$ \sigma_{gD}(T) = \{ \lambda \in \mathbb{C}: \lambda I-T   \text{ is not generalized Drazin invertible} \},$\\

$\sigma_{lgD}(T):=\{\lambda \in \mathbb{C} : \lambda I-T   \text{ is not left generalized Drazin invertible} \},$\\ 

and\\

$\sigma_{rgD}(T):=\{\lambda \in \mathbb{C} : \lambda I-T \text{ is not right generalized Drazin invertible} \}$.\\ 

It is well known that these spectra are compact sets in the complex plane, and we have,
$$\sigma_{gD} (T)=\sigma_{lgD}(T) \cup \sigma_{rgD}(T)\subset \sigma_{D}(T)=\sigma_{lD}(T) \cup \sigma_{rD}(T),$$

 $$ \sigma_{lgD}(T)\subset \sigma_{lD}(T)\subset \sigma_{ap} (T),$$
 
 and 
 $$\sigma_{rgD}(T)\subset\sigma_{rD}(T)\subset \sigma_{su} (T),$$

where

$$\sigma_{ap}(T):=\{\lambda \in \mathbb{C} : \lambda I-T  \text{ is not bounded below } \}$$
and 
$$\sigma_{su}(T):=\{\lambda \in \mathbb{C} : \lambda I-T  \text{ is not surjective} \},$$

are respectively the  approximate point spectrum  and the surjective spectrum of $T$. 

An operator  $T\in \mathcal{L}(X)$, $T$ is said to be semi-regular if $R(T)$ is closed and $N(T^{n}) \subseteq R(T)$, for all $n\in \mathbb{N}$. An  important class of operators which involves the concept of semi-regularity is the class of  operators admits a generalized Kato decomposition.
\begin{definition} (see \cite{Aiena}) $T\in \mathcal{L}(X)$  is said
to admit a generalized Kato decomposition,  abbreviated as GKD, if there exists a pair of  closed subspaces $(M,N)$ such that $(M, N) \in Red(T)$ with $T_M $ is semi-regular and $T_N$ is quasi-nilpotent. 

The pair $(M,N)$ is called the generalized Kato decomposition of   $T$ and  denoted by GKD$(M,N)$.
\end{definition}
 If we assume in the definition above that $T_N$ is nilpotent, then there exists $d\in \mathbb{N}$ for which $(T_N)^{d}=0$. In this case $T$ is said to be of Kato type operator  of  degree $d$.  
 Examples of  operators admits a generalized Kato decomposition,  are Kato  type operators, semi-regular operators, semi-Fredholm operators, quasi-Fredholm operators and generalized Drazin  invertible operators, some other examples   may be found in \cite{labro1}.

For  operator $T$  admits a generalized Kato decomposition  we have the following properties of $K(T)$ and $H_{0}(T )$.
\begin{theorem}[\cite{Aiena}] \label{thm1}
Assume that $T \in \mathcal{L} (X)$,  admits a
GKD $(M,N)$. Then 
\begin{enumerate}
	\item $K(T) = K(T_M)$ and $K(T)$ is closed.
	\item $K(T)\cap N(T)=N(T_M)$.
	\item $H_{0}(T )=H_{0}(T_M )\oplus H_{0}(T_N )=H_{0}(T_M ) \oplus N$.
\end{enumerate}
\end{theorem}
Let  $M$  be a subspace of  $X$ and let $X^*$  be the dual space  of $X$. As it is usual, $M^{\bot}=\{ x^*\in X^* : x^*(M)=0\}.$ Moreover, if $M$ and $N$ are closed linear subspaces of $X$ then $(M + N)^{\bot} = M^{\bot} \cap N^{\bot}$. The dual relation  $M^{\bot} + N^{\bot} = (M \cap N)^{\bot}$ is not always true, since $(M \cap N)^{\bot}$ is always closed but $M^{\bot} + N^{\bot}$ need not be closed. However, a classical theorem establishes
that $M^{\bot} \cap N^{\bot}$ is closed in $X^*$ if and only if  $ M + N$ is closed in $X$, (see \cite[Theorem 4.8, Chapter IV]{Kato}).
\begin{theorem}[\cite{Aiena}] \label{ttt} Let $T\in \mathcal{L}(X)$. If  $(M,N)$ is a GKD of $T$, then $(N^{\bot},M^{\bot})$ is a GKD of its adjoint $T^*$. Furthermore, if $T$ is of  a Kato type operator  then  $T^*$ is  also of  a Kato type.
\end{theorem}
For every operator $T\in \mathcal{L}(X)$, let us define the semi-regular spectrum,  the Kato  spectrum and  the generalized Kato spectrum as follows: \\

$\sigma_{se}(T):=\{\lambda \in \mathbb{C} : \lambda I-T \text{ is not semi-regular} \}$\\

$\sigma_{k}(T):=\{\lambda \in \mathbb{C} : \lambda I-T \text{ is not of Kato type} \}$\\

$\sigma_{gk}(T):=\{\lambda \in \mathbb{C} : \lambda I-T \text{ does not admit a generalized Kato decomposition} \}$\\

Recall that all the three sets defined above  are  always  compact subsets of the complex plane, (see \cite{Aiena}, \cite{Jiang}) and ordered by :
\[
\sigma_{gk}(T)\subseteq \sigma_{k}(T)\subseteq \sigma_{se}(T).
 \]
Furthermore,   the generalized Kato spectrum of an operator differs from the semi-regular spectrum on at most countably many isolated points, more precisely  the sets $ \sigma_{se}(T) \setminus \sigma_{gk}(T)$, $ \sigma_{se}(T) \setminus \sigma_{k}(T)$ and  $ \sigma_{k}(T) \setminus \sigma_{gk}(T)$ are at most countable (see \cite{Aiena} and  \cite{Jiang}).  

Note that $\sigma_{gk}(T)$  (resp. $\sigma_{k}(T)$) is not necessarily non-empty. For example, a quasinilpotent (resp. nilpotent) operator  $T$ has empty generalized Kato spectrum (resp. Kato spectrum). Furthermore, the comparison between this spectra and the spectra defined by the Drazin inverses gives

 $$\sigma_{gk}(T)\subset \sigma_{lgD}(T)\subset \sigma_{lD}(T)\subset \sigma_{ap} (T),$$
 
 and 
 $$\sigma_{gk}(T)\subset\sigma_{rgD}(T)\subset\sigma_{rD}(T)\subset \sigma_{su} (T). $$
\begin{definition}
Let  $T\in \mathcal{L}(X)$. The
operator $T$ is said to have the single-valued extension property at $\lambda_0 \in \mathbb{C}$,
abbreviated $T$ has the SVEP at $\lambda_0$, if for every neighborhood $\mathcal{ U}$ of $\lambda_0$ the
only analytic function $f :\mathcal{ U} \rightarrow X$ which satisfies the equation
$$(\lambda I- T)f(\lambda )=0$$
is the constant function $f\equiv 0$.\\
The operator $T$ is said to have the SVEP if $T$ has the SVEP at every $\lambda \in \mathbb{C}$.
\end{definition}
Trivially, an operator $T$ has the SVEP at every point of the resolvent set $\rho (T)$. Moreover, from the identity theorem for analytic functions it easily follows that $T$ has the SVEP at every point of
the boundary $\partial \sigma(T)$ of the spectrum. Hence, we have the implications:
\begin{enumerate}
\item Every operator $T$ has the SVEP at an isolated point of
the spectrum.
\item If $\lambda \notin$ acc$\sigma_{ap} (T)$, then $T$ has the SVEP at  $\lambda$.
\item If $\lambda \notin$ acc$\sigma_{su} (T)$, then  $T^*$ has the SVEP at  $\lambda$
\end{enumerate}
In particular, it has been showed that if $\lambda I-T$ admits a generalized Kato decomposition, then implications (2) and (3) may be reversed. For more properties of the SVEP, we can see \cite{LauNeu}.
\section{Left and right generalized Drazin invertible operators and the SVEP}\label{sec:3}
Now we give a characterization of the left  (resp. the right) generalized Drazin invertible  operators  in terms of generalized Kato decomposition and  the single-valued extension property. 

\begin{theorem}\label{thmKGDSVEP} An operator $T \in \mathcal{L} (X)$ is left  generalized Drazin invertible if and only if  $T$  admits a GKD $(M,N)$ and $T_M$ has the SVEP at $0$.
\end{theorem}
\begin{proof} By definition \ref{def 1.1} a left  generalized Drazin invertible operator   $T$  admits a GKD $(M,N)$ 
with $H_{0}(T)=N$ is closed, hence $T$ has the SVEP at $0$. So, $T_M$  has the SVEP at $0$.
Conversely,  if   $T$ admits a GKD$(M,N)$ with $T_M$ has the SVEP at $0$. Then by \cite[Theorem 3.14]{Aiena} $T_M$ is injective and $H_{0}(T)=N$. Since $R(T_M)$ is closed, $T_M$ is bonded below. Hence $T$ is left generalized Drazin invertible.
\end{proof}
Dually, by Definition \ref{def 1.2} and \cite[Theorem 3.15]{Aiena}, we get the following result,
\begin{theorem}\label{thmKGDSVEP1}
 An operator $T \in \mathcal{L} (X)$ is right  generalized Drazin invertible if and only if  $T$  admits a GKD $(M,N)$ and $T^*$ has the SVEP at $0$.
 \end{theorem}
Again by \cite[Theorem 3.14]{Aiena}, there are an equivalent proprieties to $T_M$ has the SVEP at $0$ for operators  admits a GKD, so we can say more about the left  generalized Drazin invertible operators. 
\begin{theorem}\label{thmKGD1} An operator $T \in \mathcal{L} (X)$ is left  generalized Drazin invertible if and only if  $T$  admits a GKD $(M,N)$ and satisfies  one of the following equivalent assertions: 
\begin{itemize}
	\item[(i)] $T$ has the SVEP at $0$,
	\item[(ii)] $T_M$ has the SVEP at $0$, 
	\item[(iii)] $T_M$ is injective,
	\item[(iv)] $H_{0}(T)=N$,
	\item[(v)] $H_{0}(T)$ is closed,
	\item[(vi)] $K(T)\cap H_{0}(T) =\{0\}$,
	\item[(vii)] $K(T)\cap H_{0}(T)$ is closed.
\end{itemize}
\end{theorem}
Similarly, by \cite[Theorem 3.15]{Aiena} we have:
\begin{theorem}\label{thmKGD2} An operator $T \in \mathcal{L} (X)$ is right  generalized Drazin invertible if and only if  $T$  admits a GKD $(M,N)$ and
satisfies  one of the following equivalent assertions: 
\begin{itemize}
	\item[(i)] $T^*$ has the SVEP at $0$,
	\item[(ii)] $T_M$ is surjective, 
	\item[(iii)] $K(T)=M$,
	\item[(iv)] $X=K(T)+ H_{0}(T)$,
	\item[(v)] $K(T)+ H_{0}(T)$ is norm dense in $X$.
\end{itemize}
\end{theorem}
The following result expresses  a characterization of the isolated points of $\sigma_{ap} (T)$ in terms of generalized Kato decomposition and  the SVEP. 
\begin{proposition}\label{prop:lgdsvep} Let $ T \in \mathcal{L} (X)$ and $0\in \sigma_{ap} (T)$. Then $0$ is an isolated point in $\sigma_{ap} (T)$ if and only if   $T$ admits a GKD$(M,N)$ and  $T_M$ has the SVEP at $0$.  
\end{proposition}
\begin{proof} Suppose that $0$ is an isolated point in $\sigma_{ap} (T)$, then $T$ has the SVEP at $0$ and by \cite[Proposition 9.]{GoMbOu}, $H_0(T)$ and $K(T)$ are closed subspaces of $X$ with $K(T)\neq X$,  $H_0(T)\neq\{0\}$ and  $K(T)\cap H_0(T)=\{0\}$. If $K(T)\oplus H_{0}(T)=X$, then $0$ is also isolated point in $\sigma (T)$ and clearly $T$ admits a GKD$(K(T),H_0(T))$. Now, assume that $K(T)\oplus H_{0}(T)\subsetneq X$, and  denote by $X_0 =K(T)\oplus H_{0}(T)$.  Observe that $X_0$ is a Banach space and  $(K(T),H_0(T))$ is a GKD of $T$ on $X_0$. So by Theorem \ref{ttt} $(H_0(T)^{\bot},K(T)^{\bot})$ is also  a GKD of $T^{*}$ over $X_{0}^{*}=H_0(T)^{\bot}\oplus K(T)^{\bot}$. In the other hand we have $H_0(T)^{\bot}+ K(T)^{\bot}=X^{*}$ and the adjoint of the inclusion map $i: X_{0} \rightarrow X$ is a  map from $X^{*}$ onto $X_{0}^{*}$ with kernel $X_{0}^{\bot}$. This implies that  $T^*$ admits a GKD over $X^*$. Again by Theorem \ref{ttt} $T$ admits a GKD$(M,N)$ viewed as a restriction  of the adjoint of $T^*$ on $X$. Further,   $T_M$ has the SVEP at $0$ because the  SVEP is inherited by the restrictions on invariant subspaces. Conversely,  if   $T$ admits a GKD$(M,N)$ and $T_M$ has the SVEP at $0$. Then by \cite[Theorem 3.14]{Aiena} $T_M$ is injective and $H_{0}=N$. Since $R(T_M)$ is closed, $T_M$ is bonded below. Hence $T$ is left generalized Drazin invertible. By  \cite[Theorem 3.8] {miloud} $0$ is an isolated point in $\sigma_{ap} (T)$.



\end{proof}
\begin{remark}If $X$ is a Hilbert space, then by \cite[Th\`eor\'eme 4.4]{Mbekhta000}  an isolated  point $\lambda$ in $\sigma_{ap} (T)$ is a singularity of the generalized resolvent, equivalently, $\lambda I-T$ admits a generalized Kato decomposition, so the   proof  above is more direct in this case.
\end{remark}
\begin{proposition}\label{prop:rgdsvep} Let $ T \in \mathcal{L} (X)$ and $0\in \sigma_{su} (T)$. Then $0$ is an isolated point in $\sigma_{su} (T)$ if and only if   $T$ admits a GKD$(M,N)$ and $T^{*}_{N^{\bot}}$ has the SVEP at $0$.
\end{proposition}
\begin{proof}Since $\sigma_{su} (T)=\sigma_{ap} (T^*)$, we apply a same analysis as in the above proof using in particular  \cite[Theorem 3.15]{Aiena} instead \cite[Theorem 3.14]{Aiena} and \cite[Theorem 3.10] {miloud} instead \cite[Theorem 3.8] {miloud}.

\end{proof}
The basic existence results of generalized Drazin inverses and  their  relation  to the SVEP, the quasinilpotent part and the analytical core are summarized in the following theorems.

\begin{theorem}\label{thm:lgd}
Assume that $T \in \mathcal{L} (X)$. The following
assertions are equivalent:

\begin{itemize}
\item[(i)] $T$ is left generalized Drazin invertible, 

\item[(ii)] 
 $T=T_{1}\oplus T_{2}$, with $T_{1}=T_{M}$ is left invertible operator and $T_{2}=T_{H_{0}(T)}$ is quasinilpotent operator,
 
\item [(iii)] $0$ is an isolated point in $\sigma_{ap} (T)$,
\item[(iv)] $T$ admits a GKD$(M,N)$ and $T_M$ has the SVEP at $0$, 
\item[(v)] $T$ admits a GKD$(M,N)$  and verified one of the  equivalent conditions of the Theorem \ref{thmKGD1},

\item [(vi)] there exists a bounded projection $P$ on $X$ such that $
TP=PT$, $T+P$ is bounded below,  $TP$  is  quasinilpotent and  $R(P)=H_{0}(T)$.
\end{itemize}
\end{theorem}

\begin{proof}The equivalence (i)$ \Longleftrightarrow $(ii) follows from \cite[Proposition 3.2]{miloud}, (iii)$ \Longleftrightarrow $(vi) follows from Proposition \eqref{prop:lgdsvep}, the implication (vi)$ \Rightarrow $(i) from \cite[Theorem 3.14]{Aiena} and (ii)$ \Rightarrow $(iii) follows from  \cite[Theorem 3.8] {miloud}. Now the equivalence (i)$ \Longleftrightarrow $(v) is proved in \cite[Theorem 3.1]{cvetkovic}.
\end{proof}
We know that the properties to be right generalized Drazin invertible or to be left generalized Drazin invertible are dual to  each other, (see \cite[Proposition 3.9]{miloud}), then we have, 
\begin{theorem}
\label{thm:rgd}
Let $T \in \mathcal{L} (X)$. The following assertions are equivalent:

\begin{itemize}
\item[(i)] $T$ is right generalized Drazin invertible,

\item[(ii)] $T=T_{1}\oplus T_{2}$, with $T_{1}=T_{K(T)}$ is right invertible
operator and $T_{2}=T_{N}$ is quasinilpotent operator,

\item [(iii)] $0$ is an isolated point in $\sigma_{su} (T)$, 

\item[(iv)] $T$ admits a GKD$(M,N)$ and $T^*$ has the SVEP at $0$,

\item[(v)] $T$ admits a GKD$(M,N)$ and satisfied one of the  equivalent conditions of the Theorem \ref{thmKGD2},

\item [(vi)] there exists a bounded projection $P$ on $X$ such that $
TP=PT$, $T+P$ is surjective,  $TP$  is  quasinilpotent and  $N(P)=K(T)$.
\end{itemize}
\end{theorem}
\begin{proof}The equivalence (i)$ \Longleftrightarrow $(ii) follows from\cite[Proposition 3.4]{miloud}, (iii)$ \Longleftrightarrow $(vi) follows from Proposition \eqref{prop:rgdsvep}, the implication (vi)$ \Rightarrow $(i) from \cite[Theorem 3.15]{Aiena} and (ii)$ \Rightarrow $(iii) follows from  \cite[Theorem 3.10] {miloud}. Now the equivalence (i)$ \Longleftrightarrow $(v) is proved in \cite[Theorem 3.2]{cvetkovic}.
\end{proof}
\begin{remark}Theorems \ref{thm:lgd} and \ref{thm:rgd} are improvement of \cite[Theorem 3.8, Theorem 3.10]{miloud}. 
\end{remark}
A direct consequence of Theorem \ref{thm:lgd}, Theorem \ref{thm:rgd}  and \cite[Lemma 3.13]{Aiena}:
\begin{corollary}Let $T \in \mathcal{L} (X)$.
\begin{itemize}
\item If $T$ is  left generalized Drazin invetible, then $H_{0} (T)^{\bot} =K(T^*).$ 
\item If $T$ is  right generalized Drazin invetible, then $K(T^*) =H_{0} (T^*)^\bot .$
\item If $T$ is  left generalized Drazin invetible, then $T^*$ is  a right generalized Drazin invetible with $N(P^*)=K(T^*) =H_{0} (T^*)^\bot=R(P)^\bot $, where $P$ is the bounded projection given in   Theorem \ref{thm:lgd}-(vi). 
\end{itemize}
\end{corollary}
Denote by 
$$ \mathcal{S}(T)=\{\lambda \in \mathbb{C}: T \text{ does not have the SVEP at } \lambda\}.$$
\begin{corollary}Let $T \in \mathcal{L} (X)$.
\begin{itemize}
\item $\sigma_{gk} (T)\cup \mathcal{S}(T)=\sigma_{lgD} (T).$
\item $\sigma_{gk} (T)\cup \mathcal{S}(T^*)=\sigma_{rgD} (T).$
\item If $T^*$ has SVEP, then $\sigma_{lgD} (T)=\sigma_{gD} (T)$ and  $\sigma_{gk} (T)=\sigma_{rgD} (T).$
\item If $T$ has SVEP, then $\sigma_{rgD} (T)=\sigma_{gD} (T)$ and  $\sigma_{gk} (T)=\sigma_{lgD} (T).$
\item If $X$ is a Hilbert space and $T$ is a self-adjoint operator,  then $\sigma_{lgD} (T)=\sigma_{rgD} (T)=\sigma_{gD} (T).$  
\end{itemize}
\end{corollary}
Similarity, for operators of Kato type we have,
\begin{corollary}Let $T \in \mathcal{L} (X)$.
\begin{itemize}
\item $\sigma_{k} (T)\cup \mathcal{S}(T)=\sigma_{lD} (T).$
\item $\sigma_{k} (T)\cup \mathcal{S}(T^*)=\sigma_{rD} (T).$
\item If $T^*$ has SVEP, then $\sigma_{lD} (T)=\sigma_{D} (T)$ and  $\sigma_{k} (T)=\sigma_{rD} (T).$
\item If $T$ has SVEP, then $\sigma_{rD} (T)=\sigma_{D} (T)$ and  $\sigma_{k} (T)=\sigma_{lD} (T).$  
\end{itemize}
\end{corollary}
In the following, we show that both $\sigma_{lgD}(T)$ and $\sigma_{rgD}(T)$ are stable under additive commuting finite rank operators. 
\begin{proposition}\label{fin}Let $T\in \mathcal{L}(X)$ and  $F$ 
is a finite rank operator on $X$ such that $TF=FT$. Then 
\begin{equation*}
\sigma_{lgD}(T+F)=\sigma_{lgD}(T).
\end{equation*}
\end{proposition}
\begin{proof}From \cite[Lemma 2.3]{lah}  we know that $acc\sigma_{ap}(T+F)=acc\sigma_{ap}(T)$. Then $\lambda \notin acc\sigma_{ap}(T+F) $ if and only if  $\lambda \notin acc\sigma_{ap}(T)$. Hence $\lambda I-(T + F)$ is left generalized Drazin invertible if and only if $\lambda I-T$ is left generalized Drazin invertible. So $\sigma_{lgD}(T+F)=\sigma_{lgD}(T).$ 
\end{proof}
As a consequence of Proposition \ref{fin} we have
\begin{proposition}Let $T\in \mathcal{L}(X)$ and  $F$ 
is a finite rank operator on $X$ such that $TF=FT$. Then 
\begin{equation*}
\sigma_{rgD}(T+F)=\sigma_{rgD}(T).
\end{equation*}
\end{proposition}

\begin{proposition}\label{propek}   Let $T\in \mathcal{L}(X)$ and $0\in \rho(T)$. Then $\lambda \in \sigma_{lgD}(T)$ if and only if $\lambda \neq 0 $ and  $ \lambda^{-1} \in \sigma_{lgD}( T^{-1})$.
\end{proposition}
\begin{proof}  we have 
$$\lambda I-T=-\lambda (\lambda^{-1}I -T^{-1})T.$$
Since $0\in \rho(T)$  and $T$ commute with $(\lambda^{-1} -T^{-1})$, it follows   that $(\lambda^{-1}I -T^{-1})_M$ is  left invertible if and only if $(\lambda I-T)_M $ is  left invertible and   $(\lambda^{-1} I-T^{-1})_N$ is quasinilpotent if and only if $(\lambda I-T)_N $ is quasi-nilpotent. This is equivalent to the statement of the Theorem.
\end{proof}
\begin{theorem} \label{thm04} Let $T,S \in \mathcal{L}(X) $. If $0\in \rho(T)\cap \rho(S)$, such that $ T^{-1}- S^{-1} $ is  finite rank operator commuting with $T$ or $S$, then
$$ \sigma_{lgD}(T)=\sigma_{lgD}(S).$$
\end{theorem}
\begin{proof}
Proposition  \ref{fin} implies that $\sigma_{lgD}(T^{-1}) = \sigma_{lgD}(S^{-1})$ , and by Proposition \ref{propek} we have $\sigma_{lgD}(T) = \sigma_{lgD}(S)$.
\end{proof}

\begin{proposition}Let $T\in \mathcal{L}(X)$ and $0 \in \rho(T)$. Then $\lambda \in \sigma_{rgD}(T)$ if and only if $ \lambda \neq 0$ and  $ \lambda^{-1} \in \sigma_{rgD}(T^{-1}).$
\end{proposition}
\begin{theorem}Let $T,S \in \mathcal{L}(X) $. If $0 \in \rho(T)\cap \rho(S)$, such that $ T^{-1}- S^{-1} $ is  a finite rank operator commuting with $T$ or $S$, then
$$ \sigma_{rgD}(T)=\sigma_{rgD}(S).$$
\end{theorem}
\begin{theorem}Let $R,T,U \in \mathcal{L}(X)$ be such that $TRT = TUT$. Then
\begin{equation*}
\sigma_{lgD}(TR)=\sigma_{lgD}(UT).
\end{equation*}
\end{theorem}
\begin{proof} Since $\sigma_{ap}(TR)\backslash \{0\}=\sigma_{ap}(UT)\backslash \{0\}$, from \cite[Theorem 1]{coor}, then it is enough to show that $TR$ is left generalized Drazin invertible  if anf only if $ UT$ is. Assume that $0 \notin \sigma_{lgD}(TR),$
then $0 \in iso\sigma_{ap}(TR)$. Therefore $TR -\mu I$ is bounded below for all small $\mu \neq 0$. Hence
$UT-\mu I$ is bounded below for all small $\mu \neq 0$. So $0 \in iso\sigma_{ap}(UT)$. Hence $UT$ is left generalized Drazin invertible if and only if $TR $ is left generalized Drazin invertible.
\end{proof}
By duality, we have:
\begin{theorem}Let $R,T,U \in \mathcal{L}(X)$ be such that $TRT = TUT$. Then
\begin{equation*}
\sigma_{rgD}(TR)=\sigma_{rgD}(UT).
\end{equation*}
\end{theorem}
In particular if $R = U$ we get
\begin{corollary}\label{co1} Let $R,T \in \mathcal{L}(X)$ then
\begin{equation*}
\sigma_{lgD}(TR)=\sigma_{lgD}(RT) \quad \text{and}~ 
\sigma_{rgD}(TR)=\sigma_{rgD}(RT).
\end{equation*}
\end{corollary}
\begin{example}Let $R, T \in \mathcal L(X)$ and $A$ be the operator defined on
$X \oplus X$ by
$$A=\left( \begin{array}{ccc}
0 & T \\
R & 0
\end{array} \right),$$
then $A^{2}=\left( \begin{array}{ccc}
TR & 0 \\
0 & RT
\end{array} \right)= TR \oplus RT$. Thus $\sigma_{lgD}(A^{2})=\sigma_{lgD}(TR) \cup \sigma_{lgD}(RT)$ which equals to $\sigma_{lgD}(TR)$ from Corollary \ref{co1}. Therefore
$\sigma_{lgD}(A)= (\sigma_{lgD}(TR))^{1/2}.$
Similarly we have $\sigma_{rgD}(A)= (\sigma_{rgD}(TR))^{1/2}.$
\end{example}
\section{Generalized Drazin inverse and local spectral theory}\label{sec:4}
We know that  if  $T\in \mathcal{L}(X)$ is  not invertible then $T$  is   generalized Drazin invertible if and only if  $X=K(T)\oplus H_0(T)$ and, with respect tho this  decomposition,
\begin{equation*}
T=T_{1}\oplus T_{2}, \text{ with } T_{1}=T_{K(T)} \text{ is invertible and } T_{2}=T_{H_0(T)} \text{ is quasi-nilpotent.}
\end{equation*}
Note that the   generalized Drazin inverse $T^D$  of  $T$, if it exists, is uniquely determined and represented, with respect of the same decomposition, as the direct sum 
\begin{equation*}
T^D=T^{-1}_{1}\oplus 0, \quad \text{ with } T^{-1}_{1} \text{ is  the inverse of  } T_{1},
\end{equation*}
Furthermore, the nonzero part of the
spectrum of $T^D$ is given by the reciprocals of the nonzero points of the spectrum of $T$, i.e.,
\begin{equation}\label{equ:opinv}
 \sigma (T^D)\setminus \{0\}=\{ \frac{1}{\lambda}: \lambda \in \sigma (T)\setminus \{0\} \}.
 \end{equation}
 Since the spectral mapping theorem holds for the approximate spectrum and the surjective spectrum, we have
 \begin{equation}
 \sigma_{ap} (T^D)\setminus \{0\}=\{ \frac{1}{\lambda}: \lambda \in \sigma_{ap} (T)\setminus \{0\} \}  
 \end{equation}
 and 
 \begin{equation}
 \sigma_{su} (T^D)\setminus \{0\}=\{ \frac{1}{\lambda}: \lambda \in \sigma_{su} (T)\setminus \{0\} \}.
 \end{equation}
 
An interesting question given in \cite{AienaT} is that there is a reciprocal relationship between the nonzero part of the local spectrum of a Drazin invertible operator  and   the nonzero part of the local spectrum of its  Drazin inverse. In the sequel we study this question in the case of the generalized Drazin invertible operators.
 
 Before this down, we shall give the relevant definitions concerning the local spectral theory. Given a bounded linear operator $T\in \mathcal{L}(X)$ , the local resolvent set $\rho_{T} (x)$ of $T$ at a point $x \in X$ is defined as the union of all open subsets $U$ of $\mathbb{C}$ such that there exists an analytic function $f : U \longrightarrow X $ satisfying
 $$ (\lambda I - T)f(\lambda)=x \quad \text{ for all } \lambda \in U.$$
 
The local spectrum $\sigma_{T} (x)$ of $T$ at $x$ is the set defined by $\sigma_{T} (x) := \mathbb{C} \setminus \rho_{T} (x)$. Obviously, $\sigma_{T} (x)\subseteq \sigma (T)$. 

The SVEP for $T$ is equivalent to saying that $\sigma_{T} (x)=\emptyset$ if and only if $x = 0$, see \cite[Proposition 1.2.16]{LauNeu}.
Note that if $T$ has SVEP then a spectral theorem holds for the local spectrum, i.e., if $f$ is an analytic function defined on an open neighborhood $U$ of $\sigma (T)$ then
$$f(\sigma_{T} (x))=\sigma_{f(T)} (x) \quad \text{ for all } x\in X.$$
See also \cite{Vrbova}.

An important invariant subspace in local spectral theory is given by the local spectral subspace of $T$ associated at a subset $\Omega \subseteq \mathbb{C}$, defined as
$$
X_{T}(\Omega)=\{x\in X: \sigma_{T} (x)\subset \Omega\}.
$$
Obviously, for every closed set $\Omega \subseteq \mathbb{C}$ we have
$$
X_{T}(\Omega)=X_{T}(\Omega\cap \sigma (T)).
$$
For a closed subset $\Omega \subseteq \mathbb{C}$, the glocal subspace $\mathcal{X}_{T}(\Omega)$ is defined as the set of all $x\in X$ for which there exists an analytic function $f: \mathbb{C}\setminus \Omega \longrightarrow X$ satisfying  $(\lambda I - T)f(\lambda)=x $ on  $\mathbb{C}\setminus \Omega$.

Obviously, for $\Omega$  a closed set, $\mathcal{X}_{T}(\Omega)\subset X_{T}(\Omega)$, and we have equality when $T$ satisfies the SVEP.

 \begin{definition} An operator $T\in \mathcal{L}(X)$ is said to have Dunford's property $(C)$, shortly property $(C)$,  if $X_{T}(\Omega)$ is closed for every closed set $\Omega \subseteq \mathbb{C}$.
 \end{definition}
 \begin{definition} An operator $T\in \mathcal{L}(X)$ is said to have  property $(Q)$,   if $H_{0}(\lambda I-T)$ is closed for every  $\lambda \in \mathbb{C}$.
 \end{definition}
Another important property which plays a central role in local spectral theory is the following one introduced by Bishop:
 \begin{definition} An operator $T\in \mathcal{L}(X)$ is said to have Bishop's property $(\beta)$, shortly property $(\beta)$,  if  for every open  set $\mathcal{U}$ of $ \mathbb{C}$ and every sequence of analytic functions $f_n : \mathcal{U} \longrightarrow X$ for which $(\lambda I-T)f_{n} (\lambda )\rightarrow 0$ uniformly on all compact subsets of  $\mathcal{U}$; then also $f_{n} (\lambda )\rightarrow 0$, again locally uniformly on $\mathcal{U}$.
 \end{definition}
 We have
 $$
\text{ property } (\beta)\Rightarrow \text{ property } (C) \Rightarrow \text{ property } (Q) \Rightarrow  \text{the SVEP }.
 $$
 See the monograph \cite{LauNeu} for a detailed study of these properties. 
 

The next first result shows that the the SVEP is  transmitted from $T$ to its  generalized Drazin inverse $T^D$,
\begin{theorem} Let $T\in \mathcal{L}(X)$ be   generalized Drazin invertible. Then $T$ has the SVEP if and only if $T^D$ has the SVEP. 
\end{theorem}
\begin{proof}Suppose that $T\in \mathcal{L}(X)$ is a generalized Drazin invertible. If   $0\notin \sigma_{} (T)$. then $f(\lambda)= \frac{1}{\lambda}$ is analytic in any open neighborhood of $\sigma (T)$ which does not contains $0$, so  by \cite[Theorem 3.3.6]{LauNeu} $T^D=T^{-1}=f(T)$ has the property SVEP. Now if $0\in \sigma (T)$. Then $T=T_{1}\oplus T_{2}$  with $T_1$ is invertible and $T_2$ is quasinilpotent.  From the first case  $T_1$ has the property SVEP, $T_1$ has also the SVEP because it is quasinilpotent. So the Drazin generalized inverse $T^D = T^{-1}_{1}\oplus 0$ has the SVEP, from \cite[Theorem 2.9]{Aiena}.

Conversely; if $T^D = T^{-1}_{1}\oplus 0$ has the SVEP then $T^{-1}_{1}$ and $T_1$ have the property SVEP. Consequently, again by \cite[Theorem 2.9]{Aiena}, $T = T_{1}\oplus T_2$ has  the property SVEP.
\end{proof}

In the following result, we show that the relation \eqref{equ:opinv} holds also in the local sens. 
\begin{theorem} Let $T\in \mathcal{L}(X)$ be   generalized Drazin invertible with generalized Drazin  inverse $T^D$. If  $T$ has the SVEP, then for every $x\in X$ we have
\begin{equation}\label{equ:locspec}
 \sigma_{T^D} (x)\setminus \{0\}=\{ \frac{1}{\lambda}: \lambda \in \sigma_{T} (x)\setminus \{0\} \}.
 \end{equation}
\end{theorem}
\begin{proof}
Suppose that $T$ has the SVEP. If  $0\notin \sigma_{} (T)$  then $T^D=T^{-1}$ and the equality \eqref{equ:locspec} follows from the spectral mapping theorem \cite[1.6]{Vrbova} applied to the function $f(\lambda)= \frac{1}{\lambda}$. Suppose that $0\in \sigma (T)$. According the
decomposition $X =K(T)\oplus H_{0}(T)$, $T_1=T_{K(T)}$ is invertible and $T_2=T_{H_{0}(T)}$ is quasi-nilpotent, then the restrictions $T_1$ and $T_2$ have the SVEP. Now, let $x \in  X$ and write $x = y + z$, with $y \in  K(T)$ and $z \in  H_{0}(T)$. Then by \cite[Theorem 2.9]{Aiena} we have
$$ \sigma_{T} (x)=\sigma_{T_1} (y)\cup \sigma_{T_2} (z).$$
The generalized Drazin  inverse $T^D = T^{-1}_{1}\oplus 0$ has the SVEP, so always  by \cite[Theorem 2.9]{Aiena}
we have
$$ \sigma_{T^D} (x)=\sigma_{T^{-1}_{1}} (y)\cup \sigma_{0} (z),$$
where
$$ \sigma_{T^{-1}_{1}} (y)=\{ \frac{1}{\lambda}: \lambda \in \sigma_{T_{1}} (y) \} \qquad \text{ for all } y\in K(T).$$
In the case $z=0$, $\sigma_{T_2} (0)=\emptyset $ and hence $\sigma_{T} (x)=\sigma_{T_1} (y)$ and, analogously, $\sigma_{T^D} (x)=\sigma_{T^{-1}_{1}} (y)$. Thus \eqref{equ:locspec}.
Now, if $z\neq 0$,  $\sigma_{T_2} (z)=\{0\}=\sigma_{0} (z)$, since  both  $T_2$ and the null  operator are  quasinilpotent operators. Furthermore,  $0\notin \sigma_{T_1} (y)$ and $0\notin \sigma_{T^{-1}_{1}} (y)$, hence  $\sigma_{T} (x)\setminus \{0\}=\sigma_{T_1} (y)$ and $ \sigma_{T^D}\setminus \{0\} (x)=\sigma_{T^{-1}_{1}} (y),$ from which we deduce,
$$ \sigma_{T^D} (x)\setminus \{0\}=\sigma_{T^{-1}_{1}} (y)=\{ \frac{1}{\lambda}: \lambda \in \sigma_{T_1} (y) \} =\{ \frac{1}{\lambda}: \lambda \in \sigma_{T} (x)\setminus \{0\} \} .$$
This complete the proof.
\end{proof}
We establish now that  that also the property $(C)$ is transferred  to the generalized Drazin inverse.
\begin{theorem}\label{thm:PC} Let $T\in \mathcal{L}(X)$ be  a generalized Drazin invertible. Then $T$ has the property $(C)$ if and only if $T^D$ the property $(C)$. 
\end{theorem}
\begin{proof}Suppose that $T\in \mathcal{L}(X)$ is  generalized Drazin invertible. If   $0\notin \sigma (T)$, then $f(\lambda)= \frac{1}{\lambda}$ is analytic in any open neighborhood of $\sigma (T)$ which does not contains $0$, so  by \cite[Theorem 3.3.6]{LauNeu} $T^D=T^{-1}=f(T)$ has the property $(C)$. Now if $0\in \sigma (T)$. Then $T$ admits a GKD $(M,N)$,  with $M=K(T)$ and $N= H_{0}(T)$,  $T_M$ is invertible and $T_N$ is quasi-nilpotent. For the Drazin generalized inverse $T^D = T^{-1}_{M}\oplus 0$, we have 
$$X_{T^D}(\Omega)=M_{T^{-1}_{M}}(\Omega)\oplus N_{0}(\Omega) \text{ for every closed set   } \Omega \subseteq \mathbb{C}.$$

Since $T_M$ is invertible, by the first case, $T_M$ has the property $(C)$ and so the inverse   $T^{-1}_{M}$ has the property $(C)$ with $M_{T^{-1}_{M}}(\Omega)$  is closed for every closed set $\Omega \subseteq \mathbb{C}$. We know that $N_{0}(\Omega) =\{0\}$ if $0\notin \Omega$ and $N_{0}(\Omega) =N$ if $0\in \Omega$. Then, $X_{T^D}(\Omega)=M_{T^{-1}_{M}}(\Omega)\oplus \{0\}$ if $0\notin \Omega$ and $X_{T^D}(\Omega)=M_{T^{-1}_{M}}(\Omega)\oplus N$ if $0\in \Omega$. In both cases $X_{T^D}(\Omega)$ is closed, and consequently $T^D$ has the property $(C)$.

Conversely; if $T^D = T^{-1}_{M}\oplus 0$ has the property $(C)$ and as above   $X_{T^D}(\Omega)=M_{T^{-1}_{M}}(\Omega)\oplus \{0\}$ if $0\notin \Omega$ and $X_{T^D}(\Omega)=M_{T^{-1}_{M}}(\Omega)\oplus N$ if $0\in \Omega$. This implies that $M_{T^{-1}_{M}}(\Omega)$ is closed and hence $T^{-1}_{M}$ has the property $(C)$. Thus $T_{M}$ has the property $(C)$. Since $X_{T}(\Omega)=M_{T_{M}}(\Omega)\oplus N_{T_N}(\Omega)$, and $T_N$  is quasinilpotent it then follows that $X_{T}(\Omega)=M_{T_{M}}(\Omega)\oplus \{0\}$ if $0\notin \Omega$, or  $X_{T}(\Omega)=M_{T_{M}}(\Omega)\oplus N$ if $0\in \Omega$. Therefore $X_{T}(\Omega)$ is closed for every closed set    $\Omega \subseteq \mathbb{C}.$ Thus $T$ has the property $(C)$.
\end{proof}
Since for operator having property $(Q)$ we have
$$H_0(\lambda I -T)=X_{T}(\{\lambda\})=\mathcal{X}_{T}(\{\lambda\}) \text{ for all  } \lambda \in \mathbb{C},$$
we can deduce from Theorem \ref{thm:PC}  that 
\begin{corollary} Let $T\in \mathcal{L}(X)$ be   generalized Drazin invertible. Then $T$ has the property $(Q)$ if and only if $T^D$ has the property $(Q)$. 
\end{corollary}
Now before to study the property $(\beta)$,  we need some preliminary results. Let $H(U,X)$ denote the space of all analytic functions from $U$ into $X$. With respect to pointwise vector space operations and the topology of locally uniform convergence, $H(U,X)$ is a Fr\'echet space. For every $T\in \mathcal{L}(X)$ and every open set $U\subseteq \mathbb{C}$, define $T_U : H(U,X) \longrightarrow H(U,X)$ by
$$(T_U f)(\lambda) := (\lambda I - T)f(\lambda) \quad \text{ for all } f \in  H(U,X) \text{ and }\lambda \in U.$$
From \cite[Proposition 3.3.5]{LauNeu}, $T$ has the property $(\beta)$ if and only if for every open set $U\subseteq \mathbb{C}$, the operator $T_U$ has closed range in $H(U,X)$. Evidently, the property $(\beta)$ is inherited by the restrictions on invariant closed subspaces. Furthermore, the  following theorem shows that  the  property $(\beta)$ is transmitted  reciprocally form a generalized Drazin invertible operator to its generalized Drazin inverse.
\begin{theorem}\label{thm:beta} Let $T\in \mathcal{L}(X)$   a generalized Drazin invertible. Then $T$ has the property $(\beta)$ if and only if $T^D$ has the property $(\beta)$. 
\end{theorem}
\begin{proof}Suppose that $T\in \mathcal{L}(X)$ is  generalized Drazin invertible and  $0\in \sigma (T)$. Then $T=T_{1}\oplus T_{2}$  with $T_1$ is invertible and $T_2$ is quasi-nilpotent. From \cite[Proposition 2.1.6]{LauNeu}, we can identify $H(U,X)$ with the direct sum $H(U,K(T))\oplus H(U,H_{0}(T))$. $T_1$  has the property $(\beta)$ and hence its inverse $T^{-1}_{1}$ has the property $(\beta)$. Now
\begin{align*}
 T^{D}_{U}[H(U,X)] &=(T^{-1}_{1}\oplus 0_U)[H(U,K(T))\oplus H(U,H_{0}(T))]\\
 &=(T^{-1}_{1})_U[H(U,K(T))] \oplus 0_U [H(U,H_{0}(T))].
\end{align*}
$$.$$
Clearly $T^{D}_{U}$ has closed range in $H(U,X)$, so $T^D$ has the property $(\beta)$.

Conversely; if $T^D = T^{-1}_{1}\oplus 0$ has the property $(\beta)$. Then as above   $T_1$ has the property $(C)$. Since the quasinilpotent operator $T_2$ has the property $(\beta)$ and the fact that

\begin{align*}
 T_{U}[H(U,X)]& =(T_{1}\oplus T_2)_U[H(U,K(T))\oplus H(U,H_{0}(T))] \\
 &=T_{1U}[H(U,K(T))] \oplus T_{2U} [H(U,H_{0}(T))];
 \end{align*} 
we conclude that  $T^{D}_{U}$ has closed range in $H(U,X)$, so $T^D$ has the property $(\beta)$.
\end{proof}
An operator $T\in \mathcal{L}(X)$ is said to have the decomposition property $(\delta)$ if the decomposition
$$X = X_T (\overline{U}) + X_T (\overline{V})$$
holds for every open cover $\{U, V \}$ of $\mathbb{C}$. Note that $T\in \mathcal{L}(X)$ has property $(\delta)$ (respectively, property $(\beta)$ ) if and only if $T^*$ has property $(\beta)$ (respectively, property $(\delta)$), see \cite[Theorem 2.5.5]{LauNeu}. If $T\in \mathcal{L}(X)$ has both property $(\beta)$ and property $(\delta)$ then $T$ is said to be decomposable.
\begin{corollary} Suppose that $T$ is generalized Drazin invertible. If $T$ has property $(\delta)$ then $T^D$ has property $(\delta)$, and analogously, if $T$ is decomposable then $T^D$ is decomposable.
\end{corollary}
\begin{proof}
 Clearly, from the definition of the generalized Drazin invertibility it follows that if $T$ is  generalized Drazin invertible then its adjoint  $T^*$ is also generalized Drazin invertible, with Drazin inverse $T^{D*}$. If $T$ has property $(\delta)$ then $T^*$ has property  $(\beta)$ and hence, by Theorem \ref{thm:beta}, also $T^{D*}$ has property $(\beta)$. By duality this implies that $T^D$ has property $(\delta)$. The
second assertion is clear: if $T$ is decomposable then $T^D$ has both properties $(\delta)$ and $(\beta)$ and the same holds for $T^D$, again by Theorem \ref{thm:beta} and the first part of the proof. Hence $T^D$  is decomposable.
\end{proof}
A natural question suggested by all the results of this section  is whether the  local spectral
properties are transmitted from a left (resp. right) generalized  Drazin invertible operator  to its left (resp. right) generalized Drazin inverse. The next example  shows that the answer to this question is negative.
\begin{example}
Let $X = \ell^2$ be the Hilbert space of all square summable complex sequences
$$x = (x_n)_n = (x_1, x_2, \ldots ), $$
indexed by the a nonnegative  integers. We define the right shift operator  $R$ and the left shift operator  $L$ in $\ell^2$  by
$$R (x_1, x_2, \ldots )=(0,x_1, x_2, \ldots ) $$
and
$$L (x_1, x_2, \ldots )=(x_2, x_3, \ldots ) .$$
We know  that  $ \sigma (R)= \sigma (L)= \mathbb{D}=\{ \lambda \in \mathbb{C}; \left|\lambda\right|\leq 1\}$ and $L=R^*$. Furthermore, $R$ is injective with colsed range and  $L$ is surjective. So $R$ is left invertible with $L$ its left inverse. Similarity,  $L$ is right invertible with $R$ its right inverse.

Now,  from  \cite[Example 1.2.8]{LauNeu}, it follows that the unilateral right shift $R$  has the property $(\beta)$ (hence has the property $(C)$, the property $(Q)$  and has the SVEP), while $L$ fails to have the SVEP, see \cite[Proposition 1.2.10]{LauNeu}.

We also have 
$$\sigma_{R} (x)= \sigma (R),$$
for every $x\in X$, so $\sigma_{R} (x)\setminus \{0\}$ is the punctured disc $\mathbb{D}\setminus \{0\}$. Consequently, the points of  $\sigma_{L} (x)\setminus \{0\}$, for any left inverse $L$, cannot be the
reciprocals of $\sigma_{R} (x)\setminus \{0\}$, otherwise $\sigma_{L} (x)$, and hence $\sigma (L)$, would be unbounded.
\end{example}
By the same notations of the definitions \ref{def 1.1} and \ref{def 1.2} and  from \cite[Proposition 1.2.10]{LauNeu} we deduce that:
\begin{proposition}
Let $T\in \mathcal{L}(X)$. We have
\begin{itemize}
\item If $T$ is  right generalized Drazin invertible and $T_{K(T)}$ has the SVEP (respectively, property $(C)$, property $(Q)$, property $(\beta)$), then $T$ is   generalized Drazin invertible.
\item If $T$ is  left generalized Drazin invertible and $T^{*}_{K(T^*)}$ has the SVEP (respectively, property $(C)$, property $(Q)$, property $(\beta)$), then $T$ is   generalized Drazin invertible.
\end{itemize}
\end{proposition}

\end{document}